\documentclass[11pt]{amsart}
\usepackage[all]{xy}
\usepackage{amsmath}
\usepackage[all]{xy}
\usepackage{amsmath}
\usepackage{amsfonts}
\usepackage{amssymb}
\usepackage{amscd}
\usepackage{amscd}
\title{Semiglobal results for $\overline\partial$ on a complex
space with arbitrary singularities}
\usepackage{amsfonts}
\usepackage{amssymb}
\newtheorem{theorem}{Theorem} [section]
\newtheorem{lemma} [theorem] {Lemma}
\newtheorem{proposition} [theorem] {Proposition}

\newtheorem{corollary} [theorem] {Corollary}

\begin{document}
\author{John Erik Forn\ae ss,  Nils \O vrelid and Sophia Vassiliadou }
\thanks{{\em 2000 Mathematics Subject Classification:} 32B10, 32J25, 32W05, 14C30}
\keywords{Cauchy-Riemann equation, Singularity, Cohomology groups}
\address{Dept. of Mathematics\\University of Michigan\\
Ann Arbor, MI 48109-1109 USA}
\address{Dept. of Mathematics\\University of Oslo\\
P.B 1053 Blindern, Oslo, N-0316 NORWAY}
\address{Dept. of Mathematics\\Georgetown University\\
Washington, DC 20057 USA} \email{fornaess@umich.edu,\;
novrelid@yahoo.no,\; sophia@math.georgetown.edu}
\begin{abstract} We obtain some $L^2$-results for the
$\overline\partial$ operator on forms that vanish to high order
near the singular set of a complex space.
\end{abstract}

\maketitle

\section{Introduction}

\medskip
\noindent Let $X$ be a pure $n-$dimensional reduced Stein space,
$A$ a lower dimensional complex analytic subset with empty
interior  containing $X_{\mbox{sing}}.$ Let $\Omega$ be an open
relatively compact Stein domain in $X$ and
$K=\widehat{\overline{\Omega}_X}$ be the holomorphic convex hull
of $\overline{\Omega}$ in $X$.  Since $X$ is Stein and
$K=\widehat{K_X},$  $K$ has a neighborhood basis of Oka-Weil
domains in $X$ (see Theorem 11, in \cite{Gunning}, Volume III,
page 102). Let $X_0$ be an Oka-Weil   neighborhood of $K$ in $X,\;
X_0\subset\subset X$. Then $X_0$ can be realized as a holomorphic
subvariety of an open polydisk in some $\mathbb{C}^N$. Set
$\Omega^*=\Omega\setminus A$. Let $d_A$ be the distance to $A$,
relative to an embedding of $X_0$ in $\mathbb C^N$ and let $|\;|$
and $dV$ denote the induced norm on $\Lambda^{\cdot} T^*\Omega^*,$
resp. the volume element (different embeddings of neighborhoods of
$\overline{\Omega}$ in $\mathbb{C}^N$ give rise to equivalent
distance functions and norms).  For a  measurable $(p,q)$  form
$u$ on $\Omega^*$  set $\|u\|^2_{N,\Omega}:=\int_{\Omega^*}|u|^2
d_A^{-N} dV$.

\medskip
\noindent In this paper we address the question of whether we can
solve the equation $\overline\partial u=f$ in $\Omega^*$ for a
$\overline\partial$-closed $(p,q)$ form $f$ on $\Omega^*$ that
vanishes to ``high order'' on $A$. Our main result is the
following theorem:

\begin{theorem}
Let $X,\,\Omega$ be as above.
For every $N_0 \geq 0,$ there exists $N \geq 0$ such that if $f$ is a
 $\overline{\partial}-$closed $(p,q)-$form on $\Omega^*,q>0,$
 with $\|f\|_{N,\Omega}<\infty,$ there is  $v \in L^{2,{\rm {loc}}}_{p,q-1}(\Omega^*)$ solving $\overline{\partial}v=f,$ with $\|v\|_{N_0,\Omega'}<\infty$
 for every $\Omega' \subset \subset \Omega$. For each $\Omega'\subset
 \subset \Omega,$ there is a solution of this kind satisfying
 $\|v\|_{N_0,\Omega'}\leq C \|f\|_{N,\Omega},$ where $C$ is a positive
constant that depends only on
 $\Omega', N,N_0.$
 \end{theorem}

\smallskip
\noindent When $A\cap \overline{\Omega}$ is a finite subset of
$\overline{\Omega}$ with \;$b\Omega\cap A=\emptyset,\;\Omega$ is
Stein and $\overline{\Omega}$ has a Stein neighborhood,  we obtain
the following corollary of Theorem 1.1:

\begin{corollary} With $N_{0},\; N$ as in Theorem 1.1 and for
$f$ a $\overline{\partial}-$closed $(p,q)-$form on $\Omega^*,q>0,$
 with $\|f\|_{N,\Omega}<\infty,$ there is a solution $u$
of $\overline\partial u=f$ on $\Omega^*$ with $\|u\|_{\Omega, N_{0}}\le
c\|f\|_{\Omega, N}$, $c$ independent of $f$. In other words, we obtain a
weighted $L^2$ estimate for $u$ on all of $\Omega$.
\end{corollary}

\smallskip
\noindent Theorem 1.1 extends to the case when $\Omega$ is just
holomorphically convex and contains a maximal compact subvariety
$B$ that is contained in $A$. It also extends to the case of
$(p,q)$ forms on $X^*$ with values in a holomorphic vector bundle
$E$ over $X$. Theorem 1.1 and Corollary 1.2 can be used to
construct analytic objects with prescribed behaviour on the
maximal, positive dimensional compact subvariety $B$ of a
holomorphically convex manifold. We also expect them to be useful
in studying the obstructions to solving $\overline\partial$ on a
deleted neighborhood of an isolated singular point of a complex
analytic set of dimension bigger than 2. The power series
arguments that were used  in the surface case in \cite{DFV},
\cite{F} might be replaced by the solution of a Cousin problem
with $L^2$ bounds which exists by our results when a finite number
of obstructions vanish.


\medskip
\noindent We have not managed to find a proof for Theorem 1.1
using  transcendental $L^2$ methods. Instead, our arguments are
based on resolution of singularities combined with cohomological
arguments in the spirit of Grauert \cite{Gr1}. In particular,
there exists a proper, holomorphic surjection $\pi:\widetilde{X}
\rightarrow X$ with the following properties:

\medskip
\noindent
i) $\widetilde{X}$ is an $n-$dimensional complex manifold.\\
ii) $\tilde{A}=\pi^{-1}(A)$ is a hypersurface in $\tilde{\Omega}$ with only ``normal
 crossing singularities", i.e. near each $x_0\in \tilde A$ there are local holomorphic coordinates $(z_1,\dots,z_n)$ in terms of which $\tilde{A}$ is given by $z_1\cdots
 z_m=0,$ where $1 \leq m \leq n$\\
iii) $\pi: \widetilde{X}\setminus \tilde{A} \rightarrow X \setminus A$ is  a biholomorphism.\\

\smallskip
\noindent This follows from the following two facts: a)\,every
reduced, complex space can be desingularized and, b) every
reduced, closed complex subspace of a complex manifold admits an
embedded desingularization. The exact statements and proofs can be
found in \cite{AHV}, \cite{BM}.

\smallskip
\noindent Let $\tilde{\Omega}:=\pi^{-1}(\Omega)$. We give
$\tilde{X}$ a Hermitian metric $\sigma$ and we consider the
corresponding distance function
$d_{\tilde{A}}(x)={\mbox{dist}}(x,\tilde{A}),$ volume element $d
\tilde{V}_{\sigma}$ and norms on $\Lambda^{\cdot} T \tilde{X}$ and
$\Lambda^{\cdot} T^*\tilde{X}.$  Let $J$ be the ideal sheaf of
$\tilde{A}$ in $\tilde{X}$ and $\Omega^p$ the sheaf of holomorphic
$p$ forms on $\tilde{X}$. We shall consider the following sheaves
on $\tilde{X}$ that are defined by:

$$
\mathcal{L}_{p,\,q}(U)=\{u\in L^{2,\,{\rm {loc}}}_{p,\,q}(U);\;\;
 \overline{\partial}u\in  L^{2,\,{\rm {loc}}}_{p,\,q+1}(U)\},
 $$

\noindent for every $U$ open subset of $\tilde{X}$ and the obvious
restriction maps $r^{U}_{V}: \mathcal{L}_{p,\,q}(U)\to
\mathcal{L}_{p,\,q}(V),$ where  $V\subset U$ are open subsets of
$\tilde{X}$. Then $u\to \overline\partial u$ defines an
$\mathcal{O}_{\tilde{X}}$-homomorphism $\overline\partial:
\mathcal{L}_{p,q}\to \mathcal{L}_{p,q+1}$ and the sequence
$$
0\to \Omega^p\to \mathcal{L}_{p,\,0}\to \mathcal{L}_{p,\,1}\to
\cdots \to \mathcal{L}_{p,\,n}\to 0
$$
\noindent is exact by the local Poincar\'e lemma for
$\overline\partial$. Since each $\mathcal{L}_{p,q}$ is closed
under multiplication by smooth cut-off functions we have a fine
resolution of $\Omega^p$. In the same way, since $J$ is locally
generated by one function, the sequence

$$
0\to J^k \Omega^p\to J^k \mathcal{L}_{p,\,0}\to \cdots \to J^k
\mathcal{L}_{p,\,n}\to 0
$$

\noindent is a fine resolution of $J^k \Omega^p$. Here, $u\in (J^k
\mathcal{L}_{p,\,q})_{x}$ if it can locally be written as $h^k
\,u_{0}$ where $h$ generates $J_x$  and $u_0\in
(\mathcal{L}_{p,\,q})_x$. It follows that

\begin{eqnarray}\label{eq: acyclicresolution}
H^q(\tilde \Omega, (J^k \Omega^p)_{|_{\Omega}})\cong \frac{
\text{ker}( \overline\partial: J^k
\mathcal{L}_{p,\,q}(\tilde\Omega)\to J^k
\mathcal{L}_{p,\,q+1}(\tilde\Omega))}{\overline\partial( J^k
\mathcal{L}_{p,\,q-1}(\tilde{\Omega}))}.
\end{eqnarray}

\medskip
\noindent Here is an outline of the proof of Theorem 1.1: The
pullback $\pi^*f$ satisfies

\begin{eqnarray}\label{eq:fakepullback}
 \int_{\tilde{\Omega}}|\pi^*f|_{\sigma}^2 d_{\tilde{A}}^{-N_1}d\tilde{V}_{\sigma}
 \leq C \int_{\Omega^*}|f|^2d_A^{-N}dV,
\end{eqnarray}

\noindent for a suitable $0<N_1<N$  and $\overline{\partial}\pi^*f=0$ on
$\tilde{\Omega}$. Suppose for
the moment that we could prove the following proposition:

\begin{proposition}
For $q>0$ and $k \geq 0$ given, there exists a natural number $\ell,\,\ell\geq k$ such that the map
$$i_*:H^q(\tilde{\Omega}, J^\ell \cdot \Omega^p)
 \rightarrow H^q(\tilde{\Omega},  J^k \cdot \Omega^p),
$$
\noindent induced by the inclusion $i: J^\ell\cdot\Omega^p\to J^k\cdot\Omega^p$, is the zero map.
\end{proposition}


\medskip
\noindent Using (\ref{eq:fakepullback}) we can show that $\pi^* f
\in J^l \mathcal{L}_{p,q}(\tilde \Omega)$ \;if \; $l\le
\frac{N_1}{2n}$. Assuming Proposition 1.3 this  means  that
$\overline\partial  v=\pi^*f$ has a solution in $J^k
\mathcal{L}_{p,\,q-1}(\tilde\Omega)$. Since $|h(x)|\le C
d_{\tilde{A}}(x)$ on compacts in the set where $h$ generates $J$
it follows that
$$
\int_{\widetilde\Omega'} |v|_{\sigma}^2 d_{\tilde{A}}^{-2k}(x)
d\tilde{V}_{\sigma}<\infty
$$

\noindent where $\widetilde{\Omega'}=\pi^{-1}(\Omega')\;\text
{and}\;\Omega'\subset\subset \Omega$.
\medskip
\noindent
Then $\overline\partial((\pi^{-1})^{*} v)=f$ on $\Omega^*$ and the final
step will be  to show that

\begin{eqnarray}\label{eq:solutionpseudoestimate}
\int_{\Omega'^{*}}  |(\pi^{-1})^{*}v|^2 d_{A}^{-N_0} dV \le c
\int_{\widetilde{\Omega'}} |v|_{\sigma}^2 d_{\tilde{A}}^{-2k}
d\tilde{V}_{\sigma}
\end{eqnarray}

\noindent when $k$ is big enough.

\medskip
\noindent {\bf{Remark:}} Proposition 1.3 was inspired  by Grauert
\cite{Gr1} (Satz 1, Section 4 ). Grauert's result corresponds to
the case where $A$ is a finite set.

\medskip
\noindent The paper is organized as follows: In section 2 we prove
Proposition 1.3. Section 3, contains the estimates for the
pullback of forms under $\pi$ and $\pi^{-1}$. In Section 4 we
prove Theorem 1.1. The proof of Corollary 1.2 is contained in
section 5. Last but not least, in section 6 we discuss some
generalizations to Theorem 1.1, Corollary 1.2.

\medskip
\noindent {\bf{Acknowledgements.}} J. E. Forn\ae ss is supported
in part by an NSF grant. Part of this work was done while the
second author was visiting the University of Michigan at Ann
Arbor, while on sabbatical leave from Oslo University. His work is
supported by the Norwegian Research Council, NFR,  and the
University of Michigan. The third author would like to thank Tom
Haines for  helpful discussions and the department of Mathematics
at the University of Oslo for its hospitality and support during
her visit in May of 2003.

\medskip
\noindent
\section{Proof of Proposition 1.3}

\noindent Following Grauert \cite{Gr1}, we consider more generally
the coherent analytic sheaves $\mathcal{S}$ on $\tilde{X}$ that
are torsion free i.e. sheaves with the property

\begin{equation}\label{eq:torsionfreeness}
T(\mathcal{S})_x=0 \;\;\;\text{for all} \;\;x\in \tilde{X}
\end{equation}

\medskip
\noindent where $T(\mathcal{S})_x=\{g_x\in \mathcal{S}_x:\;\; f_x \cdot g_x=0\;\;\text{for some}\;
f_x\neq 0,\; f_x\in \mathcal{O}_x\;\}.$

\medskip
\noindent
We shall show (Lemma 2.1) that when $\mathcal{S}$ is coherent and torsion free and $i: J^t\mathcal{S}\to \mathcal{S}$ is the inclusion  homomorphism, then the induced map $i_{\tilde{\Omega},*}: H^q(\tilde
\Omega, J^t \mathcal{S}) \to H^q(\tilde\Omega, \mathcal{S})$
is zero when $q>0$ and $t$ is big enough.
In order to exploit the idea that analytic sheaf cohomology on
$\tilde\Omega$ is concentrated over $\tilde{A}$, the exceptional set of the
resolution, we need to introduce the higher direct image sheaves, denoted by $R^q \pi_* \mathcal{S}$, of an analytic
sheaf $\mathcal{S}$ on $\tilde{X},\; q\ge 0$ and recall some basic facts about them. For $q\ge 0$ and $\mathcal{S}$
an $\mathcal{O}_{\tilde{X}}$-module, the higher direct image sheaves of $\mathcal{S}$ are the sheaves on $X$, associated to
the presheaf
$$
P: U\to H^q(\pi^{-1}(U), \mathcal{S})
$$

\noindent where $U$ open in $X$.
\medskip
\noindent
When $\phi: \mathcal{S}\to \mathcal{S}'$ is an
$\mathcal{O}_{\tilde{X}}$- homomorphism the induced maps
$\phi_{*}: H^q( \pi^{-1}(U),\mathcal{S}) \to H^q(\pi^{-1}(U),
\mathcal{S}')$, $U$ open in  $X$, determine a sheaf homomorphism
$\phi_{\#}: R^q\pi_{*}\mathcal{S}\to R^q\pi_{*} \mathcal{S}'$ on
$X$.  For future reference, we recall the $\mathcal{O}_{X}$-module
structure on $R^q\pi_*\mathcal{S}$. Given $U$  an open subset of
$X$,\, $f\in \mathcal{O}_{X}(U)$,\; we define a map $f_{U}\bullet:
\mathcal{S}_{|_{\pi^{-1}(U)}}\to \mathcal{S}_{|_{\pi^{-1}(U)}}$
described by \;$(f_{U}\bullet) s_x= (f\circ \pi)_x \cdot s_x,\;
x\in \pi^{-1}(U), \; s_x\in \mathcal{S}_x$ and let and
$(f_{U}\bullet)_{*}: H^q(\pi^{-1}(U), \mathcal{S})\to
H^q(\pi^{-1}(U), \mathcal{S})$ be the induced map on cohomology.
We can then define a map $\mathcal{O}_{X}(U)\times
H^q(\pi^{-1}(U), \mathcal{S})\to H^q(\pi^{-1}(U), \mathcal{S})$
that sends $(f,c)\in \mathcal{O}_{X}(U)\times H^q(\pi^{-1}(U),
\mathcal{S})$ to $(f_{U}\bullet)_{*} c$. It is easy to check that
it is a morphism of presheaves $\mathcal{O}_{X}(-)\times
H^q(\pi^{-1}(-), \mathcal{S})\to H^q(\pi^{-1}(-), \mathcal{S})$
which extends naturally to a  morphism on the associated sheaves
$\mathcal{O}_X\times R^q\pi_*\mathcal{S}\to R^q\pi_*\mathcal{S}$.






\medskip
\noindent
The main theorem in
Grauert \cite{Gr2}, says that the direct image sheaves $R^q\pi_*\mathcal{S}$
are coherent $\mathcal{O}_{X}$-modules, when $\mathcal{S}$ is a coherent $
\mathcal{O}_{\tilde{X}}$-module and $q\ge 0$. Since $\Omega$ is a Stein
domain, Satz 5, Section 2 in  \cite{Gr2}, gives that the natural map
$\pi_{q}: H^q(\tilde\Omega,\mathcal{S}_{|_{\tilde{\Omega}}})\to \Gamma(\Omega,
R^q\pi_*\mathcal{S}_{|_{\Omega}})$ is an isomorphism. This fact and the following lemma
will enable us to finish the proof of Proposition 1.3.

\begin{lemma}For each $q>0$ and for each coherent,
torsion free $\mathcal{O}_{\tilde{X}}$-module $\mathcal{S}$ there
exists a $t\in \mathbb{N}$ such that $i_{\tilde{\Omega},*}:
H^q(\tilde{\Omega}, J^t\mathcal{S})\to H^q(\tilde{\Omega},
\mathcal{S})$ is the zero map, where $i:
J^t\mathcal{S}\hookrightarrow \mathcal{S}$ is the inclusion map.
\end{lemma}

\noindent
\begin{proof} \noindent
We shall prove the lemma using downward induction on $q>0$.
Observe that $\tilde{\Omega}$ is an $n$-dimensional  complex
manifold with no compact $n$-dimensional connected components
since it is obtained by blow-ups from a pure $n$-dimensional Stein
space $\Omega$. It follows from the Main Theorem in Siu \cite{Siu}
that $H^{n}(\tilde{\Omega}, \mathcal{S})=0$ for every coherent
$\mathcal{O}_{\tilde{X}}$-module $\mathcal{S}$.
Hence, the statement is true for $q=n$ and any
$t\in \mathbb{N}$.


\medskip
\noindent  When $q>0,\,\,{\mbox{Supp}} R^q\pi_{*} \mathcal{S}$ is
contained in $A$. The annihilator ideal $\mathcal{A'}$ of
$R^q\pi_{*} \mathcal{S}$ is coherent and by Cartan's Theorem A
there exist functions \; $f_{1},\cdots, f_{M}\in \mathcal{A'}(X)$
that generate each stalk $\mathcal{A'}_{z}$ in a neighborhood of
$\overline{\Omega}$. Let $\mathcal{A}$ be the
$\mathcal{O}_{\tilde{X}}$-ideal generated by $\tilde f_j=f_j\circ
\pi,\;\; 1\le j\le M$.  A crucial observation which will be useful
later, is that $(\tilde{f_j})_{\tilde{\Omega},*}:
H^q(\tilde{\Omega}, \mathcal{S}_{|_{\tilde{\Omega}}})\to
H^q(\tilde{\Omega}, \mathcal{S}_{|_{\tilde{\Omega}}})$ are zero
for all $j,\,1\le j\le M, \, q>0$. To see this, consider the
following commutative diagram

\[\begin{CD}
H^q(\tilde{\Omega}, \mathcal{S}_{|_{\tilde{\Omega}}}) @>(\tilde{f_j})_{\tilde{\Omega}, *}>>   H^q(\tilde{\Omega}, \mathcal{S}_{|_{\tilde{\Omega}}})\\
@V{\cong}VV                                                                   @V{\cong}VV\\
R^q\pi_*\mathcal{S}(\Omega)@>(f_j)_{\Omega,\#}>> R^q\pi_*\mathcal{S}(\Omega) \\
\end{CD}\]

\medskip
\noindent The vertical maps are isomorphisms, due to Satz 5,
Section 2, in \cite{Gr2}. Recalling the way $\mathcal{O}_X$ acts
on $R^q\pi_*\mathcal{S}$ and using the fact that the $f_j$'s are
in the annihilator ideal of $R^q\pi_*\mathcal{S}$ we conclude that
$(f_j)_{\Omega,\#}=0$. Hence, due to the commutativity of the
above diagram $(\tilde{f_j})_{\tilde{\Omega},*}$ is zero.

\medskip
\noindent Let $Z(\mathcal{A})\,({\rm resp.}\; Z(\mathcal{A'}))$
denote the zero variety of $\mathcal{A}\,({\rm resp.}\;
\mathcal{A'})$. Since
$Z(\mathcal{A'})={\mbox{Supp}}R^q\pi_*\mathcal{S}$ is contained in
$A,$ we have that $Z(\mathcal{A})$ is contained in $\tilde{A}$
near $\overline{\tilde{\Omega}}$. Thus by R\"uckert's
Nullstellensatz for ideal sheaves, (see Theorem,  page 82 in
\cite{GrRe}), we have $J^{\mu} \subset\mathcal{A}$ on
$\tilde{\Omega}$  for some $\mu\in \mathbb{N}$.
Consider the surjection $\phi: \mathcal{S}^M \to \mathcal{A}\cdot
\mathcal{S}$ given by $(s_1,\cdots,s_M)\to \sum_{1}^{M}
\tilde{f_{j}} s_j$ and set $K=\text{ker}\phi$. Clearly, $K$  is
torsion free,  whenever  $\mathcal{S}$ is. By definition the
sequence

\begin{eqnarray}\label{eq:shortexact}
0\to K \stackrel{i}\to \mathcal{S}^{M}\stackrel{\phi}\to  \mathcal{A} \cdot \mathcal{S}\to 0
\end{eqnarray}

\noindent is exact, and it follows from (\ref{eq:torsionfreeness}) and the fact that $J$ is locally
generated by one element that
$$
0\to J^a \cdot K \stackrel{i}\to J^a \cdot \mathcal{S}^M
\stackrel{\phi}\to J^a\cdot \mathcal{A}\cdot \mathcal{S}\to 0
$$
\noindent is also exact for any $a\in \mathbb{N}$.

\medskip
\noindent
Taking all the above into consideration we obtain the following commutative diagram:



$$
\xymatrix{
 & H^q(\tilde{\Omega}, J^{a+\mu}\cdot\mathcal{S}) \ar[d] \\
 & H^q(\tilde{\Omega}, J^a\cdot \mathcal{A}\cdot
\mathcal{S})\ar[d] \ar[r]^{\delta}  & H^{q+1}(\tilde{\Omega}, J^a
\cdot K)\ar[d]_{i_1}\\
H^q(\tilde{\Omega}, \mathcal{S})^{M}\ar[dr]^{\chi}
\ar[r]^{\phi_{\tilde{\Omega}, *}}  & H^q(\tilde{\Omega},
\mathcal{A}\cdot \mathcal{S}) \ar[d]^{i_2}
\ar[r]^{\delta}  & H^{q+1}(\tilde{\Omega}, K)\\
& H^q(\tilde{\Omega}, \mathcal{S})
}
$$

\bigskip
\noindent where the third row  is exact (as part of the long exact
cohomology sequence that arises from (\ref{eq:shortexact})\,) and
the vertical maps are induced by sheaf inclusions. The map $\chi$
is defined to be $\chi:=i_2\circ \phi_{\tilde{\Omega},\,*}$ and we
can show that $\chi(c_1,\cdots,c_M)=\sum_{j=1}^{M}
(\tilde{f_j})_{\tilde{\Omega},\,*} c_j,$ where\; $c_j\in
H^q(\tilde{\Omega}, \mathcal{S}),\; 1\le j\le M$. It follows from
the induction hypothesis for $(q+1)$ applied to the coherent,
torsion-free sheaf $K$ that there exists an integer $a$ large
enough such that $i_1=0$. Then, for an element $\sigma \in
H^q(\tilde{\Omega}, \mathcal{A} \cdot \mathcal{S})$ that comes
from $H^q(\tilde{\Omega}, J^{a+\mu}\cdot \mathcal{S})$, we have
$\delta\sigma=0$, so
$\sigma=\phi_{\tilde{\Omega},*}(\sigma_1,\cdots, \sigma_M),
\;\sigma_j\in H^q(\tilde{\Omega}, \mathcal{S}),\,1\le j\le M$. By
the crucial observation above and the way $\chi$ is defined, we
conclude that $\chi$ is the zero map. Hence $i_2(\sigma)=i_2\circ
\phi_{\tilde{\Omega},*}(\sigma_1,\cdots, \sigma_m)=\sum_{j=1}^{M}
(\tilde{f_j})_{\tilde{\Omega},*} \sigma_j=0$. Thus, for
$i:J^{a+\mu}\cdot \mathcal{S}\hookrightarrow \mathcal{S}$ the
inclusion map, we have that $i_{\tilde{\Omega},*}:
H^q(\tilde{\Omega}, J^{a+\mu}\cdot \mathcal{S})\to
H^q(\tilde{\Omega}, \mathcal{S})$ is the zero map.\end{proof}

\medskip
\noindent
Choosing as $\mathcal{S}:=
J^k \Omega^p$ we obtain Proposition 1.3.

\medskip
\noindent


\medskip
\noindent
\section{Pointwise estimates for the pull back of forms under $\pi,\,\pi^{-1}$}

\medskip
\noindent Let $\sigma$ be a metric on $\tilde{X}$,
$|\;|_{x,\,\sigma}$ denote the pointwise norm of an element of
$\wedge^{r} T_{x}\tilde{X}$ or $\wedge^{r} T^{*}_{x}\tilde{X}$ for
some $r>0$ with respect to the metric $\sigma$ and $d_{\tilde{A}}$
the distance to $\tilde{A}$ in $\tilde{X}$. Let $d_{A}$ denote the
distance to $A$ relative to an embedding of a neighborhood $X_0$
of $\overline{\Omega}$ in $\mathbb{C}^N$ \;and let $|\;|_y$ denote
the pointwise norm of an element in $\wedge^{r} T_{y}(X_0\setminus
X_{\mbox{sing}})$ for some $r>0$, with respect to the restriction
of the pull back of the euclidean metric in $\mathbb{C}^N$ to
$X_0\setminus X_{\mbox{sing}}$. Let $dV,\;d\tilde{V}_{\sigma}$
denote the volume forms on $X_0 \setminus
X_{\mbox{sing}},\;\text{and}\;\tilde{X}$. The map $\pi:
\tilde{X}\setminus \tilde{A}\to X\setminus A$ is a biholomorphism
of complex manifolds. It induces a linear isomorphism $\pi_*:
\wedge^r T_x(\tilde{X}\setminus \tilde{A}) \to \wedge^r
T_{\pi(x)}(X\setminus A)$ for $x\notin \tilde{A}$.

\medskip
\noindent
\begin{lemma} We have for $x\in \widetilde{\Omega}\setminus \tilde{A},\;
v\in \wedge^{r}T_x (\widetilde{\Omega})$

\begin{eqnarray} \label{eq:pointwise}
c'\; d^{t}_{\tilde{A}}(x)&\le& d_{A}(\pi(x))\le C'\; d_{\tilde{A}}(x),\\
c\; d^{M}_{\tilde{A}}(x) \;|v|_{x,\sigma} &\le& |\pi_{*}(v)|_{\pi(x)}\le C\; |v|_{x,\sigma}.
\end{eqnarray}\\
\noindent  for some  positive constants $c',c,C',C,t,M,$ where
$c,C,M$ may depend on $r$.

\medskip
\noindent For an $r$-form $a$ in $\Omega^*$ set\;
$|\pi^{*}a|_{x,\sigma}:= {\mbox{max}} \{\;|< a_{\pi(x)},
\pi_{*}v>|\;;\; |v|_{x,\,\sigma} \le 1,\; v\in
\wedge^{r}T_{x}(\tilde{\Omega}\setminus \tilde{A})\}$, where by
$<,>$ we denote the pairing of an $r$-form with a corresponding
tangent vector. Using (7) we obtain:

\begin{equation}\label{eq:pointwiseest}
c\; d^{M}_{\tilde{A}}(x) \; |a|_{\pi(x)}\le |\pi^{*}
a|_{x,\sigma}\le C\;|a|_{\pi(x)}
\end{equation}

\noindent on $\tilde{\Omega}$, for some positive constant $M$.
\end{lemma}

\noindent {{\it Proof.}} \noindent The right hand side
inequalities in the above estimates are obvious consequences of
the differentiability of $\pi$, while the left hand side
inequalities are consequences of the Lojasiewicz inequalities (see
for example \cite{L}, or \cite{M} Chapter 4, Theorem 4.1) in the
following form:
\medskip
\noindent
\begin{lemma} Let $S$ be a real analytic subvariety of some open
subset $V$ of \;$\mathbb{R}^d$ and let $f$ be a real analytic,
real-valued function in $V$. Let $Z_f=\{x\in V; \; f(x)=0\}$.
Then, for every compact $K\subset S$, there exist positive
constants $c, m$ such that

\begin{equation*}
|f(x)|\ge c \;d(x,Z_f)^{m}
\end{equation*}
\noindent when $x\in K$.
\end{lemma}

\medskip
\noindent Lemma 3.2  generalizes easily to the case when $S$ lies
in a real analytic manifold and the distance is defined by a
Riemannian metric.

\medskip
\noindent
To prove the left hand side inequality in (6)\; let $f: \tilde{X} \times A\to
\mathbb{R}$ be given by $f(x,z)=|\pi(x)-z|^2$ and $K:= \overline{\widetilde
{\Omega}}\times\,\;(\text{compact neighborhood of} \; \overline{\Omega}\cap A$).
\medskip
\noindent
Clearly $Z_f\subset \tilde{A}\times A$. When
$x\in\overline{\widetilde{\Omega}}$ and $z$ is  the nearest point to $\pi(x)$
in $A$, we have:

$$
f(x,z)=|\pi(x)-z|^2=d(\pi(x), A)^{2}\ge c\; d((x,z), Z_f)^{m}\ge c\;
d_{\tilde{A}}(x)^{m}.
$$
\medskip
\noindent If we write $m=2t$ for some $t>0$ constant, then we
obtain from this last estimate the left hand side inequality in
(6).

\medskip
\noindent To prove the left hand side inequality in (7), we
consider the unit sphere bundle $S^{r}(\tilde{X})$ in $\wedge^{r}
T\tilde{X}$. We give $\tilde{X}$ a real analytic metric such that
$S^{r}(\tilde{X})$ becomes a real analytic manifold. We choose a
metric on $S^{r}(\tilde{X})$  such that the projection $p:
S^{r}(\tilde{X})\to \tilde{X}$ is distance decreasing. For
$\nu=(x, \xi_x)$ on the unit sphere bundle $S^{r}(\tilde{X}),$  we
set $f(\nu):=|\pi_{*} \xi_x|_{\pi(p(\nu))}^2$  and let
$K:=p^{-1}(\overline{\widetilde{\Omega}})$. Clearly, $Z_f\subset
p^{-1}(\tilde{A})$. It follows that
$|\pi_{*}\xi_x|_{\pi(p(\nu))}^{2}=f(\nu)\, \ge c
\;d(\nu,Z_f)^{L}\;\ge c\; d(p(\nu), \tilde{A})^{L}$ when $\nu\in
K$. Write $2M=L$ for some $M>0$ constant. The general case follows
by applying this last inequality to $\frac{v}{|v|_x}$ for $v\neq
0,\, v\in \wedge^r T_x( \widetilde{\Omega})$.

\medskip
\noindent Estimate (8) will be derived from (7) and the following
remark:

\medskip
\noindent {{\bf Remark:}} Let \;\;$T:V\to W$ be a linear
isomorphism of normed spaces such that $\|Tv\|\ge c \|v\|$ for
$v\in V$ and $c>0$ constant. Then $B_W(0,c)\subset T (B_V(0,1))$,
where by $B_{V}(0,1)$ we denote the unit ball in $V$ and
$B_W(0,c)$ is the ball in $W$, centered at $0$ and having radius
$c$.

\medskip
\noindent Using (7) and applying the above remark to
$\pi_*:\wedge^r T_x(\tilde{\Omega}\setminus\tilde{A})\to \wedge^r
T_{\pi(x)}(\Omega\setminus A)$ we obtain for $x\in
\tilde{\Omega}\setminus \tilde{A}, \; a_{\pi(x)}\in \wedge^r
T^*_{\pi(x)}(\Omega\setminus A)$:

\begin{eqnarray*}
|\pi^* a|_{x,\,\sigma} &=&{\mbox{max}}\{\; |< a_{\pi(x)}, \pi_{*}v>|\;\;;\;\; |v|_{x,\,\sigma} \le 1, v\in \wedge^{r}T_{x}(\tilde{\Omega}\setminus\tilde{A})\}\\
&\ge& {\mbox{max}}\{\;|<a_{\pi(x)}, w>|\;\;;\;\;|w|_{\pi(x)}\le c\, d_{\tilde{A}}(x)^M,\;w\in \wedge^{r} T_{\pi(x)}(\Omega\setminus A) \}\\
&=& c\,d_{\tilde{A}}(x)^M |a|_{\pi(x)}.\\
\end{eqnarray*}
\medskip
\noindent
This result applies in particular to the volume form in $\Omega\setminus A$  and gives:

\begin{equation}\label{eq:pullbackvolume}
c_1\, d_{\tilde{A}}(x)^{M_1} d\tilde{V}_{x,\,\sigma}\le
(\pi^*dV)_x\le C_1\, d\tilde{V}_{x,\,\sigma}.
\end{equation}

\medskip
\noindent

\section{Proof of Theorem 1.1}

\medskip
\noindent
Given $N_0\in\mathbb{N}$, choose $k\ge M+t\,\frac{N_0}{2}\ge 0$,\;with  $t,M$ as in Lemma 3.1. Then
by Proposition 1.3, there exists $\ell\ge k$ such that $H^q(\tilde{\Omega},
J^{\ell}\Omega^p)\to H^q(\tilde{\Omega}, J^{k}\Omega^p)$ is the zero homomorphism.
Choose $N\in \mathbb{N}$ such that $N\ge 2n \ell+M_1,$ where $M_1$ is as in (\ref{eq:pullbackvolume}).

\medskip
\noindent The proof of theorem 1.1 will be based on the following change of variables result:

\begin{lemma} Let $M,M'$ be orientable, Riemannian manifolds and $F:M\to M'$ an orientation
preserving diffeomorphism. Let $dV, dV'$ denote the corresponding volume elements of $M,M'$ respectively. For
$f\in L^1(M', dV')$ we have:

\begin{equation}\label{eq:changevariables}
\int_{M'} f dV'=\int_{M} (f\circ F)\; F^*(dV').
\end{equation}
\end{lemma}

\medskip
\noindent Since $\pi: \tilde{\Omega}\setminus \tilde{A}\to
\Omega\setminus A$ is a biholomorphism \& orientation-preserving
map-as long as we choose appropriate orientations on
$\Omega\setminus A,\;\;\tilde{\Omega}\setminus\tilde{A}$-, for any
$f$ satisfying $\|f\|_{N,\Omega^*}<\infty$ we have (by applying
Lemma 4.1):

$$
\int_{\Omega\setminus A} |f|^2 \,d_A^{-N} \, dV=\int_{\tilde{\Omega}\setminus \tilde{A}} |f|^2_{\pi(x)}\,
d_{A}(\pi(x))^{-N} (\pi^* dV)_x.
$$

\noindent Using the fact that

\begin{eqnarray*}
|f|_{\pi(x)} &\ge& C^{-1} |\pi^*
f|_{x,\,\sigma}\;\;\;\;\;\;({\mbox{right hand side
of (8)}}),\\
d_A (\pi(x))^{-1}&\ge&  C'^{-1}
d^{-1}_{\tilde{A}}(x)\;\;\;\;\;\;({\mbox{right
hand side of (6)),\;\; }} \\
(\pi^* dV)_{x,\,\sigma}&\ge& c_1 d^{M_1}_{\tilde{A}}(x)\,
d\tilde{V}_{x,\sigma}\;\;({\mbox{left hand side of (9)}}),
\end{eqnarray*}

\noindent we obtain

$$
\|f\|^2_{N,\Omega^*}\ge c''\; \int_{\tilde{\Omega}\setminus
\tilde{A}} |\pi^*f|^2_{x,\,\sigma}\, d_{\tilde{A}}^{{M_1}-N}
d\tilde{V}_{x,\,\sigma}
$$

\noindent for some $c''>0$ constant. Since $N$ was chosen such
that $N\ge M_1$, we see that $\overline\partial \pi^*f=0$ on
$\tilde{\Omega}$. It is not hard to show that $\pi^*f\in
J^{\ell}\mathcal{L}_{p,\,q}(\tilde{\Omega})$. By proposition 1.3
we know that there exists $v\in
J^k\mathcal{L}_{p,\,q-1}(\tilde{\Omega})$ such that
$\overline\partial v=\pi^* f$ in $\tilde{\Omega}$. Set
$u:=(\pi^{-1})^* v$. Then $\overline\partial u=f$ in $\Omega^*$
and for any $\Omega'\subset\subset \Omega$ we have:

\medskip
\noindent



\begin{eqnarray*}
\int_{\Omega'} |u|^2 d_{A}^{-N_0} dV &=&\int_{\tilde{\Omega}'\setminus \tilde{A}} |u|^2_{\pi(x)}\, d^{-N_0}_{A}(\pi(x))\, \pi^*(dV)\\
&\lesssim & \int_{\tilde{\Omega}'\setminus \tilde{A}} \;\;|v|^2_{x,\, \sigma} \;d_{\tilde{A}}^{-t N_0-2M} \, d\tilde{V}_{x,\,\sigma}\\
&\lesssim & \int_{\tilde{\Omega}'\setminus \tilde{A}}\; |v|^2_{x,\,\sigma}\; d_{\tilde{A}}^{-2k}  \, d\tilde{V}_{x,\sigma}<\infty.\\
\end{eqnarray*}

\noindent To pass from the 1st line to the 2nd one we use the fact that $|u|_{\pi(x)}\le c^{-1}d^{-M}_{\tilde{A}}(x)\,|v|_{x,\,\sigma}\;\\
d^{-N_0}_A(\pi(x))\le c'^{-N_0} d^{-tN_0}_{\tilde{A}}(x)$ and that
$(\pi^* dV)_{x, \sigma}\le C_1\, d\tilde{V}_{x,\,\sigma}$.

\bigskip
\noindent To conclude the proof of Theorem 1.1 we shall
need the following lemma:

\begin{lemma} Let $M$ be a complex manifold and let $E$ and $F$ be Frechet
spaces of differential forms (or currents) of type $(p,q-1), \;(p,q)$,
whose topologies are finer than the weak topology of currents. Assume that
for every $f\in F$, the equation $\overline\partial u=f$ has a solution
$u\in E$. Then, for every continuous seminorm $p$ on $E$, there is a
continuous seminorm $q$ on $F$ such that the equation $\overline\partial
u=f$ has a solution with $p(u)\le q(f)$ for every $f\in F,\;q(f)>0$.
\end{lemma}

\noindent
\begin{proof} Set $G=\{(u,f)\in E\times F: \overline\partial u=f\}$. Then
$G$ is closed in $E\times F$. To see this, let $(u_{\nu}, f_{\nu}) \in G$ with $u_{\nu}
\to u$ in $E$, $f_{\nu}\to f$ in $F$. For test forms $\phi\in
C^{\infty}_{0, (n-p, n-q)}(X)$ we get

\begin{eqnarray*}
\int_{M} f\wedge \phi&=&\text{lim}_{\nu\to \infty} \int_{M} f_{\nu}\wedge
\phi\\
&=&\text{lim}_{\nu\to \infty} (-1)^{p+q}  \int_{M} u_{\nu}\wedge
\overline\partial \phi=(-1)^{p+q} \int_{M} u\wedge \overline\partial\phi\\
\end{eqnarray*}

\noindent so $\overline\partial u=f$ weakly.

\medskip
\noindent
Thus, $G$ is a Frechet space and the bounded surjection $\pi_{2}: G\to F;
(u,f)\to f$ must be open. The set $\pi_{2}(\{(u,v)\in G: p(u)<1\})$ is an open
neighborhood of $0$ in $F$, and contains $\{f: q(f)\le 1\}$ for some
continuous seminorm $q$. Let $f\in F, \;0<q(f)=c$. Then $q(c^{-1}f)=1$, so
by the previous argument there exists a solution $c^{-1}u$ satisfying
$\overline\partial(c^{-1}u)=c^{-1}f$ with $p(c^{-1}u)<1$, i.e. $p(u)<c=q(f)$.\end{proof}

\medskip
\noindent
When $F$ is a Banach space with norm $\|.\|$, we conclude
that, given a seminorm $p$, there is a constant $C>0$ such that
$\{f:\;\|f\|\le C^{-1}\}\subset \overline\partial(\{u: p(u)\le 1\}$,
so $\overline\partial u=f$ has a solution $u$ with $p(u)\le C \|f\|$.
Applying this result to our situation, we see that if
$\overline\partial f=0,\; \|f\|_{\Omega,N}<\infty
$ and  $\Omega_{0}\subset \subset \Omega$, we have a solution $u$
of $\overline\partial u=f$ in $L^{2,\text{loc}}_{p,q-1}(\Omega^*)$
with $\|u\|_{\Omega_{0}, N_0}\le c \|f\|_{\Omega, N}$.

\medskip
\noindent
\section{Applications of Theorem 1.1 }

\medskip
\noindent We apply Theorem 1.1 to the case where $A\cap
\overline{\Omega}$ is a finite subset of $\overline{\Omega}$ with
\;$ b\Omega\cap A=\emptyset$, \;$\Omega\subset\subset X$ is Stein
and $\overline{\Omega} $ has a Stein neighborhood $\Omega'$.

\medskip
\noindent
\begin{proposition} With $N_{0},\; N$ as in Theorem 1.1 and $\overline\partial f=0$
on $\Omega^*$  and $\|f\|_{\Omega,N}<\infty$, there is a solution $u$
of $\overline\partial u=f$ on $\Omega^*$ with $\|u\|_{\Omega, N_{0}}\le
c\|f\|_{\Omega, N}$, $c$ independent of $f$. In other words, we obtain a
weighted $L^2$ estimate for $u$ on all of $\Omega$.
\end{proposition}

\noindent
\begin{proof} Choosing $\Omega_0\subset\subset \Omega$ containing $A\cap
\Omega$, we have a solution $u_{0}$ in $L^{2,\text{loc}}_{p,q-1}(\Omega^*)$
with $\|u_0\|_{\Omega_{0}, N_{0}}\le c \|f\|_{\Omega, N}$. We introduce a
cut-off function $\chi\in C^{\infty}(X)$ such that
$\chi=1$ on $X\setminus \Omega_{0}$ but $\chi=0$ near $A\cap \Omega$.
Set $f_{1}=\overline\partial(\chi u_{0})$. Clearly, $\|f_1\|_{L^2(\Omega)}
\le c \|f\|_{\Omega, N}$ and $f_{1}=0$ near $\Omega\cap A$.

\medskip
\noindent Let $\pi: \tilde{X}\to X$ be a desingularization of $X$
and consider the equation $\overline\partial v=\pi^{*}f_1$ on
$\tilde{\Omega }$. Let $\tilde{\Omega}_0:=\pi^{-1}(\Omega_0)$. The
equation $\overline\partial v=\pi^* f_1$ is solvable in
$L^2_{p,q-1}\,(\tilde{\Omega}_0)$. We can assume that
$\tilde{\Omega}$ can be exhausted by smoothly bounded strongly
pseudoconvex domains $\tilde{\Omega}_j:=\{z\in \tilde{\Omega};\;
\phi<c_j\;\}$ where $c_j$ are  real numbers, $\phi$ is an
exhaustion function for $\tilde{\Omega}$, of class
$C^3(\tilde{\Omega})$,  strictly plurisubharmonic outside a
compact subset and also that $b\tilde{\Omega}_0$ is smooth and
strongly pseudoconvex  and contained in each $\tilde{\Omega}_j$.
To each $\tilde{\Omega}_j$ we apply Theorem 3.4.6 in \cite{Hor}
and we obtain a solution $v_j$ to the equation $\overline\partial
v_j=\pi^*f_1$ in $\tilde{\Omega}_j$ with
$$ \int_{\tilde{\Omega}_j} |v_j|^2 \,e^{-\phi}\, d\tilde{V}_{\sigma} \le C \;\int_{\tilde{\Omega}} |\pi^* f_1|^2 \,
d\tilde{V}_{\sigma}
$$

\noindent where $C$ is a positive constant independent of $j,f$\;
(this follows from a careful inspection of the proof of Theorem
3.4.6 in \cite{Hor}).

\medskip
\noindent Consider the trivial extensions $ v_j^o$  of $v_j$
outside $\tilde{\Omega}_j$. Let $v$ be a weak limit of $v_j^o$.
Then $$\int_{\tilde{\Omega}} |v|^2 \,e^{-\phi}\,
d\tilde{V}_{\sigma}\le C \;\int_{\tilde{\Omega}} |\pi^*f_1|^2\,
d\tilde{V}_{\sigma}$$ and $\overline\partial v=\pi^* f_1$ in
$\tilde{\Omega}$. So there is a solution $v$ satisfying
$\|v\|_{L^2(\tilde {\Omega})}\le c\|f_1\|.$  Then
$w:=(\pi^{-1})^{*}v$ satisfies $\overline\partial w=f_1$ in
$\Omega^*$ but we have no longer control of its $L^2$-norm near
$A\cap\Omega$. Choose another cut-off function $\chi_{0}$ such
that $\chi_{0}=1$ on $\text{supp}\chi$ but $\chi_{0}=0$ near
$\Omega\cap A$. Then

\begin{eqnarray*}
\overline\partial( (1-\chi) u_{0}&+&\chi_{0} (\pi^{-1})^* v)=\\
&=&(1-\chi) f-\overline\partial\chi\wedge
u+\overline\partial\chi_{0}\wedge(\pi^{-1})^{*}v+\\
&+& \chi\, f+\overline\partial\chi \wedge u=\\
&=& f+\overline\partial\chi_{0}\wedge (\pi^{-1})^{*}v\\
\end{eqnarray*}

\medskip
\noindent
Finally we may solve $\overline\partial v_{1}=\overline\partial\chi_{0}
\wedge (\pi^{-1})^{*} v$ in $\Omega'^*$ (apply Theorem 1.1 to the trivial extension
of $\overline\partial\chi_{0}\wedge (\pi^{-1})^{*} v$ in $\Omega'$):

$$
\|v_{1}\|_{\Omega, N_{0}}\le c\; \|\overline\partial\chi_{0}\wedge
(\pi^{-1})^{*}v\|_{\Omega', N}\le c'\; \|\overline\partial\chi_{0}\wedge
(\pi^{-1})^{*}v\|_{L^2(\Omega)}\le C\|f\|_{\Omega, N}
$$

\noindent
since $\overline\partial\chi=0$ near $A$. Thus, $u:= (1-\chi) u_{0}+
\chi_{0} \,(\pi^{-1})^{*}v-v_{1}$ is a solution with the required
estimate. \end{proof}

\medskip
\noindent
\section{Generalizations}

\medskip
\noindent Theorem 1.1 and Corollary 1.2  also extend to the case
when $\Omega$ is a relatively compact domain in a complex space
$X$ of pure dimension $n$  with strictly pseudoconvex boundary. We
know that $\Omega$ contains a maximal positive dimensional compact
variety $B$ and let $A$ be a nowhere open  analytic subvariety of
$X$ containing $X_{{\rm sing}}$ and $B$. Then  theorem 1.1
carries over verbatim to the case described above. The proof needs
the following modifications: Let $\overline{\Omega}\subset X_{0}$
be a neighborhood with strictly pseudoconvex boundary and maximal
positive dimensional compact subvariety $B$. Take the Remmert
reduction $\phi: X_{0}\to X_{1}$ so that $X_1$ is Stein,
$\phi(B)=B_1$ is finite and $\phi: X_0\setminus B_{0} \to
X_{1}\setminus B_{1}$ is a biholomorphism. Let $\pi:
\tilde{X_{0}}\to X_0$ be a desingularization of $X_0$ such that
$\pi^{-1}(A)$ is a hypersurface with normal crossings. To obtain a
proof of Proposition 1.3 (vanishing cohomology), we need to
consider direct images $R^{q}(\phi\circ\pi)_{*} \mathcal{S}$ on
the Stein space $X_1$ and their annihilator ideal $\mathcal{A}$
for $\mathcal{S}$ coherent on $\tilde{X}$. Then, the proof carries
over.

\medskip
\noindent
Corollary 1.2, for the case when $X_{\text{sing}}\cap b\Omega$
is empty, with $A=B\cup (X_{\text{sing}}\cap \Omega)$ follows
exactly as above.

\end{document}